\documentclass[10pt,a4paper]{article}
\usepackage{amssymb, amsmath,amsfonts,amsthm}
\def\cprime{$'$}

\usepackage{graphicx}
\usepackage{color}
\usepackage[all]{xy}
\usepackage{psfrag}

\theoremstyle{plain}
\newtheorem{theorem}{Theorem}[section]
\newtheorem{proposition}[theorem]{Proposition}
\newtheorem{corollary}[theorem]{Corollary}

\newtheorem{lemma}[theorem]{Lemma}

\theoremstyle{definition}

\newtheorem{example}[theorem]{Example}

\newcommand{\R}{\mathbb{R}}
\newcommand{\Q}{\mathbb{Q}}
\newcommand{\Z}{\mathbb{Z}}
\newcommand{\N}{\mathbb{N}}
\newcommand{\UT}{\mathrm{UT}}
\newcommand{\ut}{\mathfrak{ut}}
\newcommand{\T}{T^*}
\newcommand{\g}{\mathfrak{g}}

\DeclareMathOperator{\Aut}{Aut} 
\DeclareMathOperator{\ATF}{ATF} \DeclareMathOperator{\ITF}{ITF}
\DeclareMathOperator{\GL}{GL}

\newcommand{\pf}{\textbf{Proof: }}
\renewcommand{\qed}{\hfill{$\square$}}
\newcommand{\oag}{ordered abelian group}
\newcommand{\be}{\begin{enumerate}}
\newcommand{\ee}{\end{enumerate}}
\newcommand{\beq}{\begin{equation}}
\newcommand{\eeq}{\end{equation}}
\newcommand{\ds}{\displaystyle}

\begin{document}
\bibliographystyle{plain}

\author{S. O Rourke}
\title{
Affine structures, wreath products and free affine actions on linear non-archimedean trees}\maketitle

\begin{abstract}
Let $\Lambda$ be an ordered abelian group, $\mathrm{Aut}^+(\Lambda)$ the group of order-preserving automorphisms of $\Lambda$, $G$ a group and $\alpha:G\to\mathrm{Aut}^+(\Lambda)$ a homomorphism. An $\alpha$-affine action of $G$ on a $\Lambda$-tree $X$ is one that satisfies $d(gx,gy)=\alpha_gd(x,y)$ ($x,y\in X$, $g\in G$). We consider classes of groups that admit a free, rigid, affine action in the case where $X=\Lambda$. Such groups form a much larger class than in the isometric case. We show in particular that unitriangular groups $\mathrm{UT}(n,\mathbb{R})$ and groups $T^*(n,\mathbb{R})$ of upper triangular matrices over $\mathbb{R}$ with positive diagonal entries admit free affine actions. Our proofs involve left symmetric structures on the respective Lie algebras and the associated affine structures on the groups in question. We also show that given ordered abelian groups $\Lambda_0$ and $\Lambda_1$ and an orientation-preserving affine action of $G$ on $\Lambda_0$, we obtain another such action of the wreath product $G\wr \Lambda_1$ on a suitable $\Lambda'$.

It follows that all free soluble groups, residually free groups and locally residually torsion-free nilpotent groups admit essentially free affine actions on some $\Lambda'$.

\end{abstract}
\section*{Introduction}

Let $\Lambda$ be an \oag. A $\Lambda$-metric space is defined in an entirely analogous way to a conventional ($\R$-)metric space. Just as one has $\R$-trees and group actions by isometries on $\R$-trees, group actions on $\Lambda$-trees by isometries have been studied in the literature, with a particular emphasis on free actions. This strand of geometric group theory has its origins in the notion of a length function on a group as developed by Lyndon \cite{lyndon}. Groups admitting a free isometric action (without inversions) on a $\Lambda$-tree can be thought of as a generalisation of free groups where the maximal abelian subgroups are stabilisers of lines in the $\Lambda$-tree and are naturally embedded in $\Lambda$. One of the principal novelties in contrast with free groups is that the maximal abelian subgroups are not necessarily cyclic.
See \cite{KMS-long-survey} for a survey of results in this area.

There is a more general concept of group actions by affine automorphisms on $\Lambda$-trees. An automorphism $g$ of a $\Lambda$-tree $X$ is affine with dilation factor $\alpha_g\in\Aut^+(\Lambda)$ if $d(gx,gy)=\alpha_g d(x,y)$ for all $x,y\in X$. Here $\Aut^+(\Lambda)$ denotes the group of order-preserving automorphisms of $\Lambda$. This was initiated by I. Liousse in \cite{Liousse}, who considered the case $\Lambda=\R$, and was continued by the author in \cite{affine-paper} and \cite{affine-bass} where the case of a general \oag\ was studied.

As a special case one can consider group actions on linear $\Lambda$-trees: these are $\Lambda$-trees isometric to subtrees of $\Lambda$ itself. In fact, in this paper the only $\Lambda$-trees we need to consider are of the form $X=\Lambda$ itself. Departing from standard usage, we will use the term `linear $\Lambda$-tree' in this paper to mean a $\Lambda$ tree of the form $X=\Lambda$. Any reference to an action `on $\Lambda$' should also be taken to mean an action on $\Lambda$ viewed as a $\Lambda$-tree.

In the isometric case, the analysis of groups that act on linear $\Lambda$-trees is easily done: a group admitting such an action factors through the group of isometries of $\Lambda$, which has the form $\Lambda\rtimes C_2$. It follows that groups admitting free actions on linear $\Lambda$-trees -- and more generally faithful actions -- by isometries without inversions are embeddable in $\Lambda$.

The main purpose of the current paper is to show that in the affine case groups admitting free actions on a linear $\Lambda$-tree (for some $\Lambda$) form a much richer class. Our first main result is the following.

\begin{theorem}\label{UTn-ess-free}
Every finitely generated torsion-free nilpotent group admits an essentially free
affine action on 
$\Z^m$ with the lexicographic order for some $m$.
\end{theorem}

Essential freeness of an affine action is a stronger property than freeness which implies rigidity (also called the non-nesting condition) and is robust under certain important operations, notably the base change functor and the $\Lambda$-fulfilment: see \cite[\S2.1]{affine-bass} for further details. Essentially free affine actions are automatically without inversions. Furthermore free isometric actions without inversions are essentially free.

We will use the shorthand $\ITF$ for groups that admit a free isometric action without inversions on a $\Lambda$-tree for some $\Lambda$, or $\ITF(\Lambda)$ if we wish to specify the \oag\ in question. Similarly $\ATF^e$ and $\ATF^e(\Lambda)$ will be used for groups admitting an essentially free affine action on a $\Lambda$-tree. (In \cite{affine-bass} we also use the notation $\ATF$
but we will not need this shorthand in the current paper.)

Since locally fully residually $\ATF^e$ groups are $\ATF^e$ (see \cite[Theorem 3.4]{affine-paper}), and residually torsion-free nilpotent implies fully residually torsion-free nilpotent (see \cite[\S3]{Baumslag-some-thoughts} for example), we have

\begin{corollary}\label{cor1}
Every locally residually torsion-free nilpotent group $G$ admits an essentially free
affine action on some $\Lambda$.

In
particular, right-angled Artin groups,
residually free groups and free products of residually free groups are $\ATF^e$.
\end{corollary}

A long-standing question of Baumslag asked whether free $\Q$-groups are residually torsion-free nilpotent. This has been answered recently in the affirmative by A. Jaikin-Zapirain \cite{Jaikin-Zapirain}. It follows that free $\Q$-groups have an essentially free affine action on some linear $\Lambda$.

Moreover the free action given by Corollary~\ref{cor1} is `shift-free' in the 
sense
that each dilation factor $\alpha_g$ stabilises all convex subgroups of $\Lambda$.

The interplay between residually free, fully residually free and ITF
groups is an interesting one.
Note that a free product of residually free groups is itself
residually free if and only if it is fully residually free, in which
case the constituent groups are fully residually free. Moreover a residually free group is fully residually free if and only if it contains no subgroup of the form $F\times C_\infty$ where $F$ is a free group of rank 2 (see
\cite{BBaumslag}).
Fully residually free groups are known to be $\ITF$
--- see
\cite[\S5.5]{Chiswell-book}. Moreover an ITF group is fully
residually free if and only if it is residually free. Restricting to $\Lambda=\R$, a surface group is $\ITF(\R)$ provided it is residually free. There are three `exceptional' surface groups, namely the cyclic group of order two ($\pi_1(\R\mathbb{P})$), the Klein bottle group ($\pi_1(\R\mathbb{P}\sharp\R\mathbb{P})$), and the fundamental group $$\pi_1(\R\mathbb{P}\sharp\R\mathbb{P}\sharp\R\mathbb{P})=\langle x,y,z\ |\
x^2y^2z^2=1\rangle$$ of the connected sum of three projective planes; these three groups are not residually free and do not act freely by isometries on an $\R$-tree.
On the other hand
the group $\pi_1(\R\mathbb{P}\sharp\R\mathbb{P}\sharp\R\mathbb{P})$ is $\ITF(\Z\times\Z)$.

Our next main result concerns the groups $T^*(n,\R)$ of upper triangular matrices with real entries and positive diagonal entries.

\begin{theorem}\label{Tn-ess-free}
The group $T^*(n,\R)$
admits an essentially free affine action on $\R^{k}$ for some $k$.
\end{theorem}

Our other main result concerns wreath products.

\begin{theorem}\label{wr-prods-ess-free}
Let $H$ be a group that admits an essentially free affine action on some $\Lambda_0$, and let $\Lambda_1$ be another \oag. Then the wreath product $H\wr \Lambda_1$ has an essentially free affine action on $\Lambda_0\times\Lambda'$ where $\Lambda'$ is the subgroup of the product $\prod_{\omega\in\Lambda_1}\Lambda_0$ consisting of those elements with well-ordered support.
\end{theorem}

It follows inductively that an iterated wreath product of
torsion-free abelian groups $G=A_1\wr\cdots \wr A_k$ admits a free
affine action on a linear $\Lambda$-tree. In case $k=2$ it follows
from a result of Baumslag \cite{Baumslag-wr-prods} that $G$ is
residually torsion-free nilpotent, so one can also deduce that $G$ admits a free affine action on a linear $\Lambda$-tree using Corollary~\ref{cor1}.

\begin{corollary}\label{free-poly}
Every free soluble group (i.e. every free group in the variety of
soluble groups of given derived length) admits an essentially free affine
action on a linear $\Lambda$-tree for some $\Lambda$. 
\end{corollary}

This corollary follows from Theorem~\ref{wr-prods-ess-free} together with a well-known result of {\v{S}}mel{\cprime}kin \cite{Shmelkin}.
We note also that an alternative proof can be given using Corollary~\ref{cor1}. For a result of Gruenberg asserts that
free polynilpotent groups (of given class row), and in particular free soluble groups (of given derived length), are
residually torsion-free nilpotent (see \cite[Theorem 7.1]{Gruenberg-res-props},
\cite[Theorem 1.3]{Baumslag-tf-nilpt}). Combining this result with Corollary~\ref{cor1} gives Corollary~\ref{free-poly}.
One notable difference is the action obtained via Corollary~\ref{cor1} is shift-free as already observed while the wreath product action arising from Corollary~\ref{wr-prods-ess-free}
shifts the convex subgroups of $\Lambda$ by the dilation factors $\alpha_g$.

It follows from our results that the group $\Aut^+(\Q^n)$ of order-preserving automorphisms of $\Q^n$ is $\ATF^e$ though it is not clear whether it can be shown that $\Aut^+(\R^n)$ is $\ATF^e$. The main obstacle to using the arguments in this paper in the latter case is the fact that the endomorphism ring of $\R$ is non-commutative.

In \S1 we show that the groups $\UT(n,\Z)$ are $\ATF^e(\Z^m)$ for some $m$. There is a natural action of $\UT(n+1,\Z)$ on $\Z^{n}$ by affine automorphisms. While this action is far from free, one can easily characterise those matrices 
that correspond to affine automorphisms that are rigid and have no fixed point -- we use the term essentially hyperbolic to refer to these automorphisms and matrices. The problem then reduces to showing that $\UT(n,\Z)$ admits an embedding in $\UT(m+1,\Z)$ whose image consists of essentially hyperbolic matrices. The key idea here is to consider affine structures on $\UT(n,\Q)$ starting from a left symmetric structure on $\ut(n,\Q)$. These give rise to matrices with a natural block decomposition which can then be shown are essentially hyperbolic.

We then show in \S2 how to extend this representation of $\UT(n,\Z)$ in $\UT(m+1,\Z)$ to one of $T^*(n,\R)$ in an appropriate $T^*(k,\R)$. It is necessary to increase the dimension of the codomain by $n$ and incorporate the $\log$ of the diagonal entries of $T^*(n,\R)$. 

Finally we show in \S3 we consider wreath products of groups that admit essentially free affine actions on some $\Lambda_0$ with \oag s and prove Theorem~\ref{wr-prods-ess-free}. In fact we prove a somewhat stronger result concerning a more general wreath product which we call the lexicographic wreath product which contains the restricted wreath product. 
This enables us to apply {\v{S}}mel{\cprime}kin's theorem as outlined above.

\section{Free affine actions of unitriangular groups}

\subsection{Affine actions on linear $\Lambda$-trees}
We refer to \cite{Chiswell-book} for the basic theory of
$\Lambda$-trees and isometric group actions thereon, and to
\cite[\S1]{affine-paper} for the basic theory of affine actions on
$\Lambda$-trees. Let $\Lambda$ be a linearly ordered abelian group
written additively.

Note that $\Lambda$ is itself a $\Lambda$-metric space where
$d(x,y)=|x-y|=\max\{x-y,y-x\}$.
As noted in the introduction, the only $\Lambda$-trees that we will need to consider in this paper have this form.

We will make use of the lexicographic order in a
somewhat general context. Let $\Omega$ be a linearly ordered set,
and $\Lambda_\omega$ an \oag\ for each $\omega\in\Omega$. Let
$\mathcal{L}_{\omega\in\Omega}\Lambda_\omega$ be the subgroup of the Cartesian product
$\prod_{\omega\in\Omega}\Lambda_\omega$ consisting of those
$(\lambda_\omega)_{\omega\in\Omega}$ with well-ordered support. Thus
$$\bigoplus_{\omega\in\Omega}\Lambda_\omega\leq
\underset{\omega\in\Omega}{\mathop{\mathcal{L}}}\Lambda_\omega\leq\prod_{\omega\in\Omega}\Lambda_\omega.$$ We define a partial order on the Cartesian product by declaring $(\lambda_\omega)_{\omega\in\Omega}<(\mu_\omega)_{\omega\in\Omega}$ if $\lambda_{\omega_0}<\mu_{\omega_0}$ where $\omega_0=\min\{\omega\in\Omega\ :\ \lambda_\omega\neq\mu_\omega\}$.
This restricts to a
linear order on $\mathcal{L}_{\omega\in\Omega}\Lambda_\omega$ making it an \oag, which we will
call the \emph{lexicographic product of the $\Lambda_\omega$}. If
$\Omega=\{1,\ldots,n\}$ these three groups coincide, and are
typically written $\Lambda_1\times\cdots\times\Lambda_n$, or
$\Lambda_0^n$ if the $\Lambda_\omega$ are equal to a common $\Lambda_0$. In sections 1 and 2 the indexing set $\Omega$ will be finite, though the wreath product construction in section 3 will require an infinite linearly ordered $\Omega$.

If $A$ is a subset of an \oag\ $\Lambda$ such that $x,z\in A$ and $x\leq y\leq z$ imply $y\in A$, we say that $A$ is a convex subset; in particular one may speak of a convex subgroup. The convex subgroups of an \oag\ are linearly ordered by inclusion. In the case of a lexicographic product as above, each $\omega$ determines a convex subgroup. If $x,y\in\Lambda$ and $nx<|y|$ for all $n\in\Z$ we write $x\ll |y|$; this is equivalent to the proper inclusion $[x]\subset[y]$ where $[x]$ denotes the convex subgroup spanned by $x$.

The group of $o$-automorphisms (that is, order-preserving automorphisms) of an \oag\ $\Lambda$
will be denoted $\Aut^+(\Lambda)$. Let $G$ be a group and
$\alpha:G\to\Aut^+(\Lambda)$ ($g\mapsto\alpha_g$) a homomorphism. An
\emph{$\alpha$-affine action} of $G$ on a $\Lambda$-tree $X$ is an action
on $X$ such that $d(gx,gy)=\alpha_g d(x,y)$ for all $x,y\in X$; of course isometric actions correspond to the case where $\alpha$ is trivial. 
In the case $\Lambda=\R^n$ the $o$-automorphisms of $\Lambda$ include coordinate-wise rescaling maps (by a positive dilation factor) as well as functions of the form $\mathbf{x}\mapsto\mathbf{x}+\zeta(\mathbf{x})$ where $\zeta:\R^n\to\R^n$ is an endomorphism with $\zeta(\mathbf{y})\ll |\mathbf{y}|$ for all $\mathbf{y}\neq 0$. It follows that the group $T^*(n,\R)$ of upper triangular matrices with positive diagonal entries embeds in $\Aut^+(\R^n)$ with the natural action on $\R^n$ given by left multiplication. In this context, it is useful to view elements of $\R^n$ as column vectors $(x_n,\ldots,x_1)^\top$, which we will do throughout this paper. (Here $\mathbf{x}^\top$ denotes the transpose.) Thus, for example, $\left(\begin{array}{c}0\\ \vdots\\ 0\\ 1\end{array}\right)\ll \left(\begin{array}{c}0\\ \vdots\\ 1\\ 0\end{array}\right)\ll\left(\begin{array}{c}1\\ \vdots\\ 0\\ 0\end{array}\right)$. This convention enables $o$-automorphisms as above to be represented by upper triangular matrices acting on the left.
We note however that not all $o$-automorphisms of $\R^n$ arising as described from endomorphisms can be represented as elements of $T^*(n,\R)$. See \cite[\S1.1]{affine-paper} for more details.

In this section will define a free affine action of $\UT(n,\Q)$, the group of rational upper triangular matrices with diagonal entries equal to 1,  on $\Q^m$ (for a suitable $m$) viewed as a $\Q^m$-tree.
Building on this, in section 2 we will describe a free affine action of the group $T^*(n,\R)$ on a suitable $\R^k$.

Just as the affine group can be represented as a subgroup of $\GL(n+1,\R)$, if $\alpha_g\in T^*(n,\R)$ we can represent an $\alpha_g$-affine automorphism $g:\mathbf{x}\mapsto\alpha_g\mathbf{x}+\nu_g$ of $\R^n$ in matrix form as follows. $$\left(\begin{array}{c}g\mathbf{x}\\ 1\end{array}\right)=\left(\begin{array}{c|c}\alpha_g & \nu_g\\ \hline 0 & 1\end{array}\right)\left(\begin{array}{c}\mathbf{x}\\ 1\end{array}\right)$$
Thus $\UT(n+1,\R)$ has a natural affine action on $\R^n$, as does the subgroup $\T(n,\R)\ltimes\R^n$ of $T^*(n+1,\R)$ consisting of matrices with 1 in the bottom right corner.

Note further that the map $\sigma:\mathbf{x}\mapsto-\mathbf{x}$ is an
isometry (generating a cyclic group of order 2), so we can extend
the action of $\T(n,\R)\ltimes \R^n$ above to one of
$\langle \sigma\rangle\ltimes\T(n,\R)\ltimes \R^n$.

Henceforth we will assume that all affine actions on linear $\Lambda$-trees preserve the order,
and we will identify elements of $\UT(n+1,\R)$ with the corresponding affine automorphism of $\R^n$ as described above.

An affine automorphism $g $ of $\Lambda$ is
\emph{rigid} if no subset of $\Lambda$ of the form $[\lambda_1,\lambda_2]=[\lambda_2,\lambda_1]=\{\mu\in\Lambda\ :\ \lambda_1\leq\mu\leq\lambda_2\}$ ($\lambda_1\leq\lambda_2$) is mapped properly into itself by $g ^{\pm 1}$. If $g $ preserves the orientation of $\Lambda$, as we will generally assume, this is equivalent to requiring that if $g ^{\epsilon}\lambda>\lambda$ for some $\lambda$ and $\epsilon=\pm 1$ then $g ^{\epsilon}\lambda'>\lambda'$ for all $\lambda'\in\Lambda$.

An orientation-preserving automorphism $g $ of $\Lambda$ is \emph{hyperbolic} if no point is fixed by
$g $. We note that affine hyperbolic automorphisms on more general $\Lambda$-trees are not as well-behaved as one might expect from the isometric case; see \cite[\S1]{affine-paper} for further discussion.
Here
we will consider affine actions on $\Lambda$-trees that are essentially free; see \cite[\S3.1]{affine-bass} for the general definition and basic properties of essentially free actions. In our case where the $\Lambda$-tree is equal to $\Lambda$ an affine action is essentially free if each non-trivial group element $g\in\T(n+1,\R)$ is rigid and hyperbolic (or \emph{essentially hyperbolic});
this amounts to the requirement 
that the lowest non-zero entry of $g-1$
appears in the last column, and this entry is strictly lower than
every other non-zero entry of $g-1$. (Here $1$ denotes the identity matrix.) That is, a non-trivial automorphism $g=(a_{ij})\in\T(n,\R)$ is
essentially hyperbolic if the following implication holds. \begin{center}If, for
some $i<n+1$, either $a_{ii}\neq 1$ or $a_{ij}\neq 0$ for some
$i<j<n+1$,\\ then $a_{k{(n+1)}}\neq 0$ for some
$i<k<n+1$.\end{center} Note that a product of
essentially hyperbolic automorphisms is not necessarily essentially hyperbolic.

The promised free affine action of
$G=\UT(n,\Q)$ on $\Q^{m}$ for a suitable $m$, will be obtained via an embedding
$\bar{\gamma}:G\to\UT(m+1,\Q)$ and showing that 
non-trivial elements of the image are
essentially hyperbolic with respect to the natural affine action on $\Q^{m}$. In fact our arguments in this section and the next work with $\Q$ replaced by any ordered field, with one notable exception: while the embedding of $\T(n,\R)$ in $\Aut^+(\R^n)$ is not surjective, the corresponding embedding of $\T(n,\Q)$ in $\Aut^+(\Q^n)$ is an isomorphism. Things are simpler still in case $\Lambda=\Z^n$; one has \begin{equation}\label{aff-Zn}\Aut^+(\Z^n)\ltimes\Z^n\cong\UT(n+1,\Z).\end{equation}
This will enable us to conclude that every group admitting an essentially free affine action on $\Z^n$ is finitely generated torsion-free nilpotent.

\subsection{Affine structures and left symmetric structures}
We will proceed by considering the standard left symmetric
structure on the Lie algebra $\mathfrak{g}=\ut(n,\Q)$ consisting of all $n\times n$ matrices with rational entries and zeroes on and below the main diagonal, and using it to obtain a complete affine structure $d\bar{\gamma}$ on $\mathfrak{g}$, and thence the promised
free action of $\UT(n,\Q)$. We refer to \cite[\S2]{Dekimpe-Malfait} for an account of the requisite Lie algebra material (including an example covering the case $n=3$).
Essentially we follow the steps (1)-(4)
as described in that paper in reverse order. We note however that not all simply transitive actions arising in this way are essentially hyperbolic in our sense: the final example of \cite{Dekimpe-decat} is not essentially hyperbolic.

So observe first that the Lie algebra
$\mathfrak{g}=\mathfrak{ut}(n,\Q)$ is graded:
that is, $\g$ decomposes as a direct sum of $n-1$ vector subspaces
$\sigma_i$ consisting of matrices $S_i^*$
whose non-zero entries appear only on the $i$th superdiagonal.
Crucially, one has $[\sigma_i,\sigma_j]\subseteq \sigma_{i+j}$
where, as usual, $[A,B]=AB-BA$, and $\sigma_k$ is taken to be the
zero subalgebra if $k\geq n$. In fact in our case we have
$\sigma_i\sigma_j\subseteq \sigma_{i+j}$.

Arising from any grading, one can define a left symmetric structure
that is compatible with the Lie bracket using the following formula and extending the operation to $\mathfrak{g}$ by bilinearity.
\begin{equation}\label{lsa}S_i^*\cdot
S_j^*=\frac{j}{i+j}[S_i^*,S_j^*]\end{equation} 
Next, define a linear map
$\lambda':\g\to\mathfrak{gl(g)}$ via
\begin{equation}\lambda'(\mathbf{x}):\mathbf{y}\mapsto
\mathbf{x}\cdot\mathbf{y}.\end{equation}

In order to represent the
$\lambda'(\mathbf{x})$ as matrices we now represent $\g$ as a set of
column vectors as follows. For $\mathbf{x},\mathbf{y}\in\g$, we write $\mathbf{x}=(x_{ij})$ and $\mathbf{y}=(y_{ij})$ with $x_{ij}=y_{ij}=0$ for $i\geq j$.
We now use the following notation for the superdiagonals of $\mathbf{x}$ and $\mathbf{y}$, \begin{equation}\label{xy}\begin{array}{cc}\begin{array}{lcl}S_1&=&(x_{12},x_{23},\ldots,x_{n-1\ n})\\
S_2&=&(x_{13},x_{24},\ldots,x_{n-2\ n})\\
S_3&=&(x_{14},x_{25},\ldots,x_{n-3\ n})\\ &\vdots\\
S_{n-1}&=&(x_{1n})\end{array}
&
\begin{array}{lcl}T_1&=&(y_{12},y_{23},\ldots,y_{n-1\ n})\\
T_2&=&(y_{13},y_{24},\ldots,y_{n-2\ n})\\
T_3&=&(y_{14},y_{25},\ldots,y_{n-3\ n})\\ &\vdots\\
T_{n-1}&=&(y_{1n}),\end{array}\end{array}\end{equation} and write $\mathbf{x}$ 
in the form  $$\mathbf{x}=\prec
S_1,S_2,\ldots,S_{n-1}\succ.$$

In the opposite direction we associate an element $S^*$ of $\g$ to a vector $S$ of length $n-i$ via $$S^*=\prec 0,\ldots,0,\underbrace{S}_{\begin{array}{c}i\mbox{th}\\ \mbox{position}\end{array}},0,\ldots,0\succ.$$

Then $S_i^*\cdot T_j^*=(W_{i,j})^*$ where
$W_{i,j}=(w_1,\ldots,w_{n-(i+j)})$ and \begin{equation}\label{w-ij}w_k=\frac{j}{i+j}\left(x_{k,i+k}y_{i+k,i+j+k}-y_{k,j+k}x_{j+k,i+j+k}\right).\end{equation}
Of course if $n\leq i+j$, we take $(W_{i,j})^*$ to be the zero matrix.

We now set $\lambda'(\mathbf{x})(\mathbf{y})=\mathbf{x}\cdot\mathbf{y}$, and writing $\mathbf{x}$ and $\mathbf{y}$ as linear combinations of $S_i^*$ and $T_j^*$ respectively ($1\leq i,j\leq n-1$), we obtain \label{x-dot-y}\begin{eqnarray*}\lambda'(\mathbf{x})(\mathbf{y})
&=&\sum_{r=2}^{n-1}\sum_{i=1}^{r-1}
(W_{i,r-i})^*\\ &=& \sum_{r=2}^{n-1}
\bar{W}_r^*\end{eqnarray*} where
$\bar{W}_r=(\bar{w}_1,\ldots,\bar{w}_{n-r})$ and
\begin{equation}\label{wk}\bar{w}_k
=\sum_{i=1}^{r-1}\frac{r-i}{r}\left(x_{k,i+k}y_{i+k,r+k}-y_{k,r-i+k}x_{r-i+k,r+k}\right).\end{equation}
Note that all entries of $\bar{W}_r^*$ are zero apart from the
entries on the $r$th superdiagonal.

Taking $m=\frac{n(n-1)}{2}$ we define a linear isomorphism
$t:\g\to\Q^{m}$ via \begin{equation}\label{t-iso}\prec S_1,S_2,\ldots,S_{n-1}\succ\
\mapsto(S_{n-1},S_{n-2},\ldots,S_1)^\top.\end{equation} Now $\lambda'(\mathbf{x})$
can be represented as an $m\times m$ matrix. More precisely,
$\lambda'(\mathbf{x})$ determines a matrix
$\lambda(\mathbf{x})\in\mathfrak{ut}(m,\Q)$ where
$\lambda(\mathbf{x})=t\cdot \lambda'(\mathbf{x})\cdot  t^{-1}$.

Thus
$$\begin{array}{rcll}\lambda(\mathbf{x})(T_{n-1},\ldots,T_1)^\top&=&\lambda(\mathbf{x})\cdot
t\prec T_1,\ldots,T_{n-1}\succ & \\ &=&
t\cdot \lambda'(\mathbf{x})\prec T_1,\ldots,T_{n-1}\succ & \\ &=&
t\cdot \lambda'(\mathbf{x})(\mathbf{y}) & \\
&=& t \displaystyle{\sum_{r=2}^{n-1} \bar{W}_{r}^*} &
\\ &=&
\\ &=& t\displaystyle{\sum_{r=2}^{n-1}}\prec 0,\ldots, 0,\underbrace{\overline{W_{r}}}_r,0,\ldots,0\succ
\\ &=& \displaystyle{\sum_{r=2}^{n-1}}(0,\ldots, 0,\underbrace{\overline{W_{r}}}_{n-r},0,\ldots,0)^\top\\ &
\\ &=& (\overline{W_{n-1}},\overline{W_{n-2}},\ldots,\overline{W_{2}},0)^\top.
\end{array}$$

Now observe that $T_{n-j}$ is a vector of length $j$ consisting of entries of the form $y_{k,k+(n-j)}$. Moreover we can write $\lambda(\mathbf{x})$ in the block form \begin{equation}\label{Mblock}\lambda(\mathbf{x})=\left(\begin{array}{ccccccc}
M_{n-1,n-1} & M_{n-1,n-2} & \cdot & \cdot & \cdot & M_{n-1,1}\\ M_{n-2,n-1} & M_{n-2,n-2} & \cdot & \cdot &
\cdot & M_{n-2,1} \\ \cdot & \cdot &
\cdot & \cdot & \cdot & \cdot \\ \cdot & \cdot & \cdot
& \cdot & \cdot & \cdot \\ \cdot & \cdot & \cdot & \cdot &
\cdot & \cdot \\ M_{2,n-1} & M_{2,n-2} & \cdot & \cdot
& \cdot & M_{2,1} \\ M_{1,n-1} & M_{1,n-2} &
\cdot & \cdot & \cdot & M_{1,1}
\end{array}\right)\end{equation}

where each $M_{n-i,n-j}$ is an $i\times j$ matrix.
It follows that \begin{equation}\label{barW}\sum_{j=1}^{n-1}M_{n-i,n-j}T_{n-j}^\top=\bar{W}_{n-i}^\top.\end{equation}

In the next section we will need a somewhat detailed description of the entries of $\lambda(\mathbf{x})$, so we will derive it here. First, for $\nu\in\N$  put $k_\nu=\min\{k\in\N: \nu\leq\sum_{\eta=1}^k \eta\}$, and $r_\nu=\nu-\sum_{\eta=1}^{k_\nu-1}\eta$.  Thus $1\leq r_\nu\leq k_\nu$ and \begin{equation}\label{block-index}\begin{array}{rl}\rho &=\sum_{\eta=1}^{k_\rho-1}\eta+r_\rho\\
\sigma &=\sum_{\eta=1}^{k_\sigma-1}\eta+r_\sigma.\end{array}\end{equation}

\begin{lemma}\label{m-rho-sigma}The $(\rho,\sigma)$ entry of $\lambda(\mathbf{x})$ is \begin{equation}m_{\rho,\sigma}=\left\{\begin{array}{cl}\ds{\frac{n-k_\sigma}{n-k_\rho}x_{r_\rho,r_\rho+(k_\sigma-k_\rho)}} & \mbox{if }k_\rho<k_\sigma\mbox{ and }r_\rho<r_\sigma \mbox{ and }r_\sigma-r_\rho=k_\sigma-k_\rho\\ \\ \ds{-\frac{n-k_\sigma}{n-k_\rho}x_{n-k_\sigma+r_\rho,n-k_\rho+r_\rho}} & \mbox{if } k_\rho<k_\sigma\mbox{ and }r_\rho=r_\sigma\\ 0 & \mbox{otherwise.}\end{array}\right.\end{equation}
\end{lemma}

\pf
Observe that the $(\rho,\sigma)$ entry of $\lambda(\textbf{x})$  is the $(r_\rho,r_\sigma)$ entry of the $(k_\rho, k_\sigma)$ block $M_{n-k_\rho,n-k_\sigma}$ of $M$: recall that this block has dimensions  $k_\rho\times k_\sigma$. Using equation (\ref{barW})
we now obtain \begin{equation}\label{M-S}
\sum_{j=1}^{n-1}M_{n-k_\rho,n-j}\cdot \left(\begin{array}{c}y_{1,n-j+1}\\ y_{2,n-j+2}\\ \vdots\\ y_{j,n}
\end{array}\right)=\bar{W}_{n-k_\rho}^\top=
\left(\begin{array}{c}z_1\\ z_2\\ \vdots\\ z_{k_\rho}
\end{array}\right)
\end{equation}
where by~(\ref{wk}) \begin{equation}\label{z-p}
z_p=\sum_{l=1}^{n-k_\rho-1}\frac{n-k_\rho-l}{n-k_\rho}
\left(x_{p,p+l}y_{p+l,n-k_\rho+p}-y_{p,n-k_\rho+p-l}x_{n-k_\rho+p-l,n-k_\rho+p}\right).
\end{equation}

Suppose that $T_{n-k_\sigma}$ is the only non-zero superdiagonal of $\mathbf{y}$ and that all entries of this superdiagonal are non-zero. Then (\ref{z-p}) simplifies to $M_{n-k_\rho,n-k_\sigma}T_{n-k_\sigma}^\top=\bar{W}_{n-k_\rho}^\top$. Suppose that the $(r_\rho,r_\sigma)$ entry of $M_{n-k_\rho,n-k_\sigma}$ is non-zero. Then equating the $r_\rho$ row of~(\ref{M-S}) and the expression~(\ref{z-p}) with $p=r_\rho$ we see that the coefficient of $y_{r_\sigma,r_\sigma+(n-k_\sigma)}$ is non-zero whence
$n-k_\rho-l=n-k_\sigma$ giving $l=k_\sigma-k_\rho$. Moreover we have either \be\item $r_\rho<r_\sigma$ and the coefficient of $y_{r_\sigma,r_\sigma+(n-k_\sigma)}\neq 0$ has the form $\ds{\frac{n-k_\sigma}{n-k_\rho}x_{r_\rho,r_\rho+(k_\sigma-k_\rho)}}$ or
\item $r_\rho=r_\sigma$ and the coefficient of $y_{r_\sigma,r_\sigma+(n-k_\sigma)}\neq 0$ has the form $\ds{-\frac{n-k_\sigma}{n-k_\rho}x_{n-k_\sigma+r_\rho,n-k_\rho+r_\rho}}$.
    \ee

It now follows that $m_{\rho,\sigma}$ has the claimed form.\qed

We now obtain a
function $d\bar{\gamma}:\mathfrak{g}\to\mathfrak{ut}(m+1,\R)$
via
\begin{equation}\label{dbargamma}d\bar{\gamma}(\mathbf{x})=\left(\begin{array}{cc}\lambda(\mathbf{x}) &
t(\mathbf{x})\\ 0 & 0\end{array} \right).
\end{equation}

Note that $d\bar{\gamma}$ is a Lie algebra homomorphism into the affine algebra $\mathfrak{aff}(m)$ where elements of the image have linear part $\lambda(\mathbf{x})$ and translation part $t(\mathbf{x})$, where $\lambda(\mathbf{x})$ is a nilpotent matrix for all $\mathbf{x}$ and where $t$ is an isomorphism of vector spaces. That $d\bar{\gamma}$ preserves the Lie bracket follows from its definition in terms of the left symmetric structure. (It does not, however, preserve matrix multiplication.) It follows that $\bar{\gamma}=\exp \cdot d\bar{\gamma}\cdot\log$ is a (Lie) group homomorphism $\UT(n,\Q)\to\UT(m+1,\Q)$;
here $\exp$ and $\log$ denote the matrix exponential and logarithm functions,
which can be evaluated as a finite sum in the case at hand since
$d\bar{\gamma}(\mathbf{x})$ is a nilpotent matrix.

\begin{lemma}\label{admissible}
\be\item
If a matrix of the form $d\bar{\gamma}(\mathbf{x})$ is non-zero,
then at least one block of its last column is non-zero. Moreover, if
$S_{i_0}\neq 0$ and $S_{i}=0$ for all $i<i_0$ then $M_{i,j}=0$
for all $i$ and $j$ such that $i-j<i_0$.
\item $\bar{\gamma}(g)$
is essentially hyperbolic for $g\neq 1$.\ee
\end{lemma}

\pf
(i) If all blocks of the last column of $d\gamma(\mathbf{x})$ are zero then $\mathbf{x}=0$ since $t$ is a linear isomorphism and $d\gamma(\mathbf{x})=0$. The stated condition on $S_{i_0}$ implies that it is the $i_0$th superdiagonal and the first non-zero superdiagonal of $\mathbf{x}$.
Since any non-zero entries of $M_{n-i,n-j}$ are constant multiples of entries of the $(j-i)$th superdiagonal of $\mathbf{x}$, we see that the blocks $M_{i,j}$ are zero for $i-j< i_0$.

(ii)
For $g\neq 1$ we put $\mathbf{x}=\log g$, take $i_0$ as in part (1), and note that the foregoing argument shows that all block superdiagonals (with respect to the block decomposition in (\ref{Mblock})) of $N=d\bar{\gamma}(\mathbf{x})$ below the $i_0$th are zero. Thus $N^k$ is a matrix in which all block superdiagonals below the $(ki_0)$th are zero. It follows that the lowest non-zero superdiagonal of $\sum_{k=0}^{m-1}\frac{1}{k!}N^k-I$ remains the $i_0$th, and the lowest non-zero entry in this matrix is in the last column and is strictly lower than any other non-zero entries.

Since
$\bar{\gamma}(g)=
\exp(N)$, this shows that $\bar{\gamma}(g)$ is essentially hyperbolic.\qed

\begin{example}\label{n=4}
Let us illustrate the constructions described by considering the case $n=4$. Here $m=6$. Consider generic elements $\mathbf{x}$ and $\mathbf{y}$ of $\mathfrak{ut}(4,\Q)$; these have the form $$\mathbf{x}=\left(\begin{array}{rrrr}0& x_{12} & x_{13} & x_{14}\\ 0 & 0 & x_{23} & x_{24}\\ 0 & 0 & 0 & x_{34}\\ 0 & 0 & 0 & 0\end{array}\right)\mbox{ and }\mathbf{y}=\left(\begin{array}{rrrr}0& y_{12} & y_{13} & y_{14}\\ 0 & 0 & y_{23} & y_{24}\\ 0 & 0 & 0 & y_{34}\\ 0 & 0 & 0 & 0\end{array}\right).$$ Expressing $\mathbf{x}$ and $\mathbf{y}$ as sums of superdiagonal matrices, applying the definition of the left symmetric structure~(\ref{lsa}) and using the bilinearity of this structure we obtain \begin{equation}\lambda'(\mathbf{x})(\mathbf{y})=\mathbf{x}\cdot\mathbf{y}=\left(\begin{array}{cccc}0& 0 & \frac{1}{2}(x_{12}y_{23}-y_{12}x_{23}) & \frac{2}{3}(x_{12}y_{24}-y_{13}x_{34})+\frac{1}{3}(x_{13}y_{34}-y_{12}x_{24}) \\ 0 & 0 & 0 & \frac{1}{2}(x_{23}y_{34}-y_{23}x_{34}) \\ 0 & 0 & 0 & 0\\ 0 & 0 & 0 & 0\end{array}\right).
\end{equation}
Take $t:\mathfrak{ut}(4,\Q)\to\Q^6$ to be the map $\mathbf{y}\mapsto(y_{14}, y_{13}, y_{24}, y_{12}, y_{23}, y_{34})^\top$. Then the map $\lambda'(\mathbf{x})$ can be represented as the matrix $$\lambda(\mathbf{x})=\left(\begin{array}{c|cc|ccc}0 & -\frac{2}{3}x_{34} & \frac{2}{3}x_{12} & -\frac{1}{3}x_{24} & 0 & \frac{1}{3}x_{13}\\ \hline 0 & 0 & 0 & -\frac{1}{2}x_{23} & \frac{1}{2}x_{12} & 0 \\ & & & & & \\ 0 & 0 & 0 & 0 & -\frac{1}{2}x_{34} & \frac{1}{2}x_{23}\\ \hline 0 & 0 & 0 & 0 & 0 & 0\\
0 & 0 & 0 & 0 & 0 & 0\\ 0 & 0 & 0 & 0 & 0 & 0\end{array}\right)\in
\mathfrak{ut}(6,\Q).$$
One can easily verify that the formula in Lemma~\ref{m-rho-sigma} is satisfied. The matrix $d\bar{\gamma}(\mathbf{x})$ is obtained by adjoining the column $t(\mathbf{x})=(x_{14}, x_{13}, x_{24}, x_{12}, x_{23}, x_{34})^\top$ and a row of zeroes.

Recall that the exponential of an element $N$ of $\mathfrak{ut}(n,\Q)$ has the form $\sum_{k=0}^{n-1}\frac{1}{k!}N^k$, while in a similar spirit $\log(I+B)$ may be similarly evaluated as $\sum_{k=1}^{n-1}\frac{(-1)^{k+1}}{k} B^{k}$ for $B\in\ut(n,\Q)$. 

Therefore for $A=\left(\begin{array}{cccc}1 & c & e & f\\ 0 & 1 & b & d\\ 0 & 0 & 1 & a\\ 0 & 0 & 0 & 1\end{array}\right)$, we have 
$$\begin{array}{rl}\bar{\gamma}(A)&=\exp \cdot d\bar{\gamma} \cdot \log(A)=\exp \cdot d\bar{\gamma}\left(\begin{array}{cccc}0 & c & e-\frac{1}{2}bc & f-\frac{1}{2}cd-\frac{1}{2}ae+\frac{1}{3}abc\\ 0 & 0 & b & d-\frac{1}{2}ab\\ 0 & 0 & 0 & a\\ 0 & 0 & 0 & 0\end{array}\right)\\ \\
\end{array}$$

$$\begin{array}{rl}
&=\exp\left(\begin{array}{c|cc|ccc|c}0 & -\frac{2}{3}a & \frac{2}{3}c & -\frac{1}{3}d+\frac{1}{6}ab & 0 & \frac{1}{3}e-\frac{1}{6}bc & f-\frac{1}{2}cd-\frac{1}{2}ae+\frac{1}{3}abc\\ \hline 0 & 0 & 0 & -\frac{1}{2}b & \frac{1}{2}c & 0 & e-\frac{1}{2}bc\\ & & & & & \\ 0 & 0 & 0 & 0 & -\frac{1}{2}a & \frac{1}{2}b & d-\frac{1}{2}ab\\ \hline 0 & 0 & 0 & 0 & 0 & 0 & c\\
0 & 0 & 0 & 0 & 0 & 0 & b\\
0 & 0 & 0 & 0 & 0 & 0 & a\\ \hline 0 & 0 & 0 & 0 & 0 & 0 & 0\end{array}\right)
\\ \\ &=\left(\begin{array}{c|cc|ccc|c}1 &
-\frac{2}{3}a & \frac{2}{3}c & -\frac{1}{3}d+\frac{1}{3}ab & -\frac{1}{3}ac & \frac{1}{3}e
& f-\frac{1}{3}cd-\frac{2}{3}ae+\frac{1}{3}abc\\ \hline 0 & 1 & 0 & -\frac{1}{2}b & \frac{1}{2}c & 0 & e-\frac{1}{2}bc\\ & & & & & \\ 0 & 0 & 1 & 0 & -\frac{1}{2}a & \frac{1}{2}b & d-\frac{1}{2}ab\\ \hline 0 & 0 & 0 & 1 & 0 & 0 & c\\
0 & 0 & 0 & 0 & 1 & 0 & b\\
0 & 0 & 0 & 0 & 0 & 1 & a\\ \hline 0 & 0 & 0 & 0 & 0 & 0 & 1\end{array}\right).
\end{array}$$

\end{example}

Note that in Example~\ref{n=4} the entries of $\bar{\gamma}(A)$ are not typically integers.
We next show that $\bar{\gamma}$ can be adjusted so that $\bar{\gamma}$ maps $\UT(n,\Z)$ into $\UT(m+1,\Z)$.

\begin{lemma}\label{embeds in UTZ}
Let $\Gamma$ be a finitely generated subgroup of $G=\UT(n,\Q)$. Then there exists a diagonal matrix $P$ such that $P\Gamma P^{-1}$
is a subgroup of $\UT(n,\Z)<\UT(n,\Q)$. In
particular, $\Gamma$ embeds in $\UT(n,\Z)$. If the elements of $\Gamma$ are essentially hyperbolic, so are those of $P\Gamma P^{-1}$.
\end{lemma}

\pf For $\kappa\neq 0$ let $D_i(\kappa)$ be the diagonal matrix with $i$th diagonal entry equal to $\kappa$ and other diagonal entries equal to 1.
Conjugating a matrix $A$
by $D_i(\kappa)$
has the effect of multiplying the $i$th row of $A$ by $\kappa$ and 
dividing the $i$th column by $\kappa$. Therefore, given a finite set of matrices $A\in \UT(n,\Q)$, such as an inverse-closed finite generating set of $\Gamma$, we can take $d_i$ to be a common multiple of all denominators of entries in the $i$th rows of the matrices $A$. Put $\bar{d_i}=\prod_{j=i}^n d_j$, and $\bar{D}_i=D_i(\bar{d_i})$. Then conjugating the matrices $A$ by $P=\bar{D}_n\bar{D}_{n-1}\cdots \bar{D}_1$ yields a set of upper triangular matrices with integer entries and each diagonal entry equal to 1. It is easy to see that essential hyperbolicity is preserved by this conjugation. 
\qed

\begin{theorem}
The following conditions are equivalent.
\be\item $G$ is finitely generated torsion-free nilpotent.
\item $G$ embeds in $\UT(n,\Z)$ for some $n$.
\item $G$ admits an essentially free affine action on a $\Z^n$-tree for some $n$.
\item $G$ admits a faithful, orientation-preserving, affine action on a $\Z^n$-tree for some $n$.
\ee
\end{theorem}

\pf
Lemmas~\ref{admissible} and \ref{embeds in UTZ} show that $\UT(n,\Z)$ admits an embedding in $\UT(m+1,\Z)$ such that non-trivial elements of the image are essentially hyperbolic.
The implication (ii)$\Rightarrow$(iii) follows.
The implication (iv)$\Rightarrow$(ii) follows from 
equation~(\ref{aff-Zn}), and  the implication (i) $\Rightarrow$ (ii) is a well-known result; see \cite[Theorem 7.5]{Edmonton}. The remaining implications are clear.\qed

\section{Free affine actions of upper triangular groups}

In this section we exhibit free rigid affine actions of the group $T^*(n,\R)$ on $\R^k$-trees where $T^*$ denotes the respective group of upper triangular matrices with positive diagonal entries. 
Our main task is to extend the actions of
the unitriangular groups $U_n=\UT(n,\R)$ described in the previous section to the groups $T^*(n,\R)$.

Note that embeddings of $\UT(n,\R)$ in $\UT(m+1,\R)$ do not necessarily extend to embeddings of the respective upper triangular groups in general, a fact brought to my attention by Yves de Cornulier and Florian Eisele in response to a question of mine on MathOverflow -- see \cite{YCor} and \cite{Eisele}.

Taking $\varphi$ to be the 
embedding of $\UT(n,\R)$ in $\UT(m+1,\R)$ of the previous section, and keeping $m=n(n-1)/2$, let us write $\varphi(u)=\left(\begin{array}{cc}\varphi_0(u)& b(u)\\ 0 & 1\end{array}\right)$ where $\varphi_0(u)\in\UT(m,\R)$ and $b(u)\in\R^m$.

Next, let $d$ be a diagonal matrix with diagonal entries $d_1,\ldots,d_n$ respectively and consider the conjugation map 
$\chi_d:x\mapsto dxd^{-1}$. Now $t(\mathbf{x})\mapsto t(\chi_d\cdot\mathbf{x})$ is a vector space automorphism of $\R^m$, which can therefore be represented as a matrix $d^*$. Direct calculation shows that $\chi_d$ has the effect of replacing the entry $x_{ij}$ of $\mathbf{x}$ by $\frac{d_i}{d_j}x_{ij}$, and therefore $d^*$ can be seen to be a block diagonal matrix whose
$i$th block is the $i\times i$ diagonal matrix with entries $d_1/d_{n-i+1}$, $d_2/d_{n-i+2}$,\ldots $d_i/d_n$ ($1\leq i\leq n-1$).

We will find it useful to denote the column vector $$\left(\begin{array}{c}\log d_1\\ \vdots\\ \log d_n\end{array}\right)$$ by $\log d$.

We now define $\bar{\varphi}:T^*(n,\R)\to T^*(m+n+1,\R)$ by the assignments
$$\begin{array}{rl}\bar{\varphi}(u) &=\left(\begin{array}{cccc}\varphi_0(u) & 0 & b(u)\\ 0 & I_n & 0 \\ 0 & 0 & 1\end{array}\right) \\ \\
\bar{\varphi}(d) &=\left(\begin{array}{cccc} d^* & 0 & 0\\ 0 & I_n & \log d\\ 0 & 0 & 1\end{array}\right)\end{array}$$
for $u\in U_n$ and $d\in D_n$ and setting $\bar{\varphi}(ud)=\bar{\varphi}(u)\bar{\varphi}(d)$
We next show that $\bar{T}^*=\langle \bar{\varphi}(U_n),\bar{\varphi}(D_n)\rangle$ is a group isomorphic to $T^*(n,\R)$ whose natural action on $\R^{m+n}$ is affine and essentially free.

For $d\in D_n$ recall that $d^*\in D_m$ as defined above. We will also write 
$\tilde{d}$ for $\left(\begin{array}{cc}d^* & 0\\ 0 & 1\end{array}\right)\in D_{m+1}$, and
$\bar{d}$ for $\bar{\varphi}(d)\in D_{m+n+1}$.

\begin{lemma}\label{10-lemma}
\be
\item\label{exp-chi} $\chi_{\tilde{d}}\cdot\exp =\exp\cdot\chi_{\tilde{d}}$.
\item\label{log-chi} $\chi_d\cdot\log =\log\cdot\chi_d$.
\item\label{t-chi-d} $d^*\cdot t=t\cdot\chi_d$
\item\label{lambda-chi} $\chi_{d^*}\cdot\lambda=\lambda\cdot\chi_d$.
\item\label{dgamma-chi} $\chi_{\tilde{d}}\cdot d\bar{\gamma}=d\bar{\gamma}\cdot\chi_d$.
\item\label{phi-chi} $\chi_{\tilde{d}}\cdot{\varphi}={\varphi}\cdot\chi_d$.
\item\label{phi-0-chi} $\chi_{d^*}\cdot\varphi_0={\varphi_0}\cdot\chi_d$;
\item\label{phi-0-chi2} $d^*\cdot b=b\cdot\chi_d$
\item\label{phi-bar-chi} $\chi_{\bar{d}}\cdot\bar{\varphi}=\bar{\varphi}\cdot\chi_d$.
\item\label{phi-bar-rest} 
$\T(n,\R)$ is isomorphic to $\bar{T}^*$
\ee
\end{lemma}

\pf
\ref{exp-chi} and \ref{log-chi} are standard (and easily checked) facts, and \ref{t-chi-d} is immediate from the definition of $d^*$.

\ref{lambda-chi}. 
Consider the $(\rho,\sigma)$ entry of $\lambda\cdot\chi_d(\mathbf{x})=\lambda(d\mathbf{x}d^{-1})$.
By Lemma~\ref{m-rho-sigma} 
$(\rho,\sigma)$ entry of $\lambda\cdot\chi_d\mathbf{x}$ is obtained by multiplying the corresponding entry of $\lambda(\mathbf{x})$ by
\begin{eqnarray*}\left\{\begin{array}{cl}\ds{\frac{d_{r_\rho}}{d_{r_\rho+(k_\sigma-k_\rho)}}} & \mbox{if }k_\sigma>k_\rho\mbox{ and }r_\sigma>r_\rho \mbox{ and }k_\sigma-k_\rho=r_\sigma-r_\rho\\ \\ \ds{\frac{d_{n-k_\sigma+r_\rho}}{d_{n-k_\rho+r_{\rho}}}} & \mbox{if } k_\sigma>k_\rho\mbox{ and }r_\sigma=r_\rho
\\ 1 & \mbox{otherwise.}
\end{array}\right.\end{eqnarray*}

On the other hand, conjugation of $\lambda(\mathbf{x})$ by $d^*$ has the effect of multiplying the $(\rho,\sigma)$ entry by $\ds{\frac{d_{r_\rho}}{d_{r_\rho+n-k_\rho}}\frac{d_{r_\sigma+n-k_\sigma}}{d_{r_\sigma}}}$. It is now clear that in each of the cases distinguished above, the $(\rho,\sigma)$ entries of $\lambda(d\mathbf{x}d^{-1})$ and of $d^*\lambda(\mathbf{x})(d^*)^{-1}$ agree. That is, $\lambda\cdot\chi_d=\chi_{d^*}\cdot\lambda$, as claimed.

\ref{dgamma-chi} follows from \ref{t-chi-d} and \ref{lambda-chi} applied to (\ref{dbargamma}).

\ref{phi-chi} follows from \ref{exp-chi}, \ref{log-chi} and \ref{dgamma-chi} applied to $\varphi=\bar{\gamma}=\exp\cdot d\bar{\gamma}\cdot\log$.

\ref{phi-0-chi}, \ref{phi-0-chi2}. We have on the one hand \begin{eqnarray*}\varphi\cdot\chi_d(u)&=&\left(\begin{array}{c|c}\varphi_0(dud^{-1}) & b(dud^{-1})\\ \hline 0 & 1\end{array}\right),\end{eqnarray*}
and on the other hand \begin{eqnarray*}\chi_{\tilde{d}}\cdot\varphi(u)&=&\left(\begin{array}{c|c}d^* & 0 \\ \hline 0 & 1\end{array}\right)\left(\begin{array}{c|c}\varphi_0(u) & b(u) \\ \hline 0 & 1\end{array}\right)\left(\begin{array}{c|c}(d^*)^{-1} & 0 \\ \hline 0 & 1\end{array}\right)\\ &=& \left(\begin{array}{c|c}d^*\varphi_0(u)(d^*)^{-1} & d^*b(u) \\ \hline 0 & 1\end{array}\right).
\end{eqnarray*}

From \ref{phi-chi} these matrices are equal, giving $\varphi_0\cdot\chi_{d}=\chi_{d^*}\cdot\varphi_0$ and $b\cdot\chi_d=d^*b$.

\ref{phi-bar-chi} This now follows from \ref{phi-0-chi} and \ref{phi-0-chi2} applied to the definition of $\bar{\varphi}$.

\ref{phi-bar-rest} To show that $\bar{\varphi}$ is a monomorphism we note first that the restriction of $\bar{\varphi}$ to $U_n$ is injective since $\varphi$ is injective, and the restriction of $\bar{\varphi}$ to $D_n$ is injective since $d\mapsto\log d$ is; both restrictions are clearly homomorphisms. It is clear that these restrictions have trivial intersection, and \ref{phi-bar-chi} implies that $\bar{\varphi}(U_n)$ is normalised by  $\bar{\varphi}(D_n)$. Moreover, by \ref{phi-bar-chi}, we have $\bar{\varphi}(dud^{-1})=\bar{\varphi}(d)\bar{\varphi}(u)\bar{\varphi}(d)^{-1}$. The result follows.
\qed

We can now deduce

\begin{theorem}
Let $n\in\N$ and $m=n(n-1)/2$.
The group $\T(n,\R)$ admits an essentially free affine action on $\R^{m+n}$ via $g\cdot \mathbf{r}=\mathbf{s}$ where $\left(\begin{array}{c}\mathbf{s}\\ 1\end{array}\right)=\bar{\varphi}(g)\left(\begin{array}{c}\mathbf{r}\\ 1\end{array}\right)$. 
\end{theorem}

\pf By Lemma~\ref{10-lemma}\ref{phi-bar-rest} $\T(n,\R)$ is isomorphic to $\bar{T}^*$. Moreover the bottom right entry of each $\bar{g}\in \bar{T}^*$ is equal to 1, so there is a natural affine action of $\bar{T}^*$ on $\R^{m+n}$ as described in \S1.1. Each $\bar{g}$ has the form $\bar{\varphi}(u)\bar{\varphi}(d)=\left(\begin{array}{ccc}d^*\varphi_0(u) & 0 & d^*b(u)\\ 0 & I_n & \log d\\ 0 & 0 & 1\end{array}\right)$. If $d\neq 1$ then $\log d\neq 0$, while if $d=1$ then $\bar{g}=\bar{\varphi}(u)$. Provided $\bar{g}$ is non-trivial it is seen to be essentially hyperbolic in either case. Therefore the action of $\bar{T}^*$, and with it the action of $\T(n,\R)$ as given in the theorem, are essentially free.
\qed

\section{Actions of wreath products}

\begin{proposition}\label{cartesian-product}
Let $\Omega$ be a linearly ordered set, and for $\omega\in\Omega$ let $\Lambda_\omega$ be an \oag\ on which $H_\omega$ admits an order-preserving $\alpha^{(\omega)}$-affine action. Set $\Lambda$ equal to the lexicographic product $\mathcal{L}_{\omega\in\Omega}\Lambda_\omega$ and $H=\mathcal{L}_{\omega\in\Omega}H_\omega$ and for $h=(h_\omega)_{\omega\in\Omega}\in H$
and $\lambda=(\lambda_\omega)_{\omega\in\Omega}$ let $\alpha_h\lambda=(\alpha_{h_\omega}^{(\omega)}\lambda_\omega)_{\omega\in\Omega}$.

The action of $H$ on $\Lambda$ given by $(h_\omega)_{\omega\in\Omega}\cdot(\lambda_\omega)_{\omega\in\Omega}=(h_\omega\cdot\lambda_\omega)_{\omega\in\Omega}$ is $\alpha$-affine and order-preserving. If the given actions of $H_\omega$ are
respectively
\be\item free;
\item rigid
\item 
essentially free
\ee 
so is that of $H$.
\end{proposition}
(See \S 1.1 for the definition of the lexicographic product.)

\pf
We will outline the straightforward proof. One checks that $\alpha$ is a homomorphism by showing that $\alpha_g(\alpha_h(\lambda_\omega)_{\omega\in\Omega})=(\alpha_{g_\omega}^{(\omega)}\alpha_{h_\omega}^{(\omega)}
\lambda_\omega)_{\omega\in\Omega}=\alpha_{gh}(\lambda_\omega)_{\omega\in\Omega}$ where $g=(g_\omega)_{\omega\in\Omega}$. 
Similarly $g\cdot h(\lambda_\omega)_{\omega\in\Omega}=(g_\omega h_\omega\lambda_\omega)_{\omega\in\Omega}=gh(\lambda_\omega)_{\omega\in\Omega}$. 

To see that the action is $\alpha$-affine, consider $\lambda'\geq\lambda''\in\Lambda$, and note that $d(h(\lambda'_\omega)_{\omega\in\Omega},h(\lambda''_\omega)_{\omega\in\Omega})
=h(\lambda'_\omega)_{\omega\in\Omega}-h(\lambda''_\omega)_{\omega\in\Omega}
=(\alpha_{h_\omega}^{(\omega)}(\lambda'_\omega-\lambda''_\omega))_{\omega\in\Omega}
=\alpha_h((\lambda'_\omega)_{\omega\in\Omega}-(\lambda''_\omega)_{\omega\in\Omega})$.

It is easy to check that if the action of each $H_\omega$ is free and orientation-preserving, then so is the action of $H$. Finally if $h(\lambda_\omega)_{\omega\in\Omega}>(\lambda_\omega)_{\omega\in\Omega}$, then there exists $\omega_0\in\Omega$ such that $h_\omega\lambda_\omega=\lambda$ for all $\omega<\omega_0$ and $h_{\omega_0}\lambda_{\omega_0}>\lambda_{\omega_0}$. If $(\lambda'_\omega)_{\omega\in\Omega}$ is another element of $\Lambda$ then $h_\omega\lambda'_\omega=\lambda'_\omega$ for all $\omega<\omega_0$, since the actions of $H_\omega$ are assumed to be rigid. Similarly $h_{\omega_0}\lambda'_{\omega_0}>\lambda'_{\omega_0}$. Thus $h(\lambda'_\omega)_{\omega\in\Omega}>(\lambda'_\omega)_{\omega\in\Omega}$.\qed\\

For a group $H$ and an \oag\ $\Lambda_1$, we denote by $H\wr_{\omega}\Lambda_1$ the subgroup of the unrestricted wreath product $H\ \bar{\wr}\ \Lambda_1$ consisting of those $(\lambda^*,(h_\lambda)_{\lambda\in\Lambda_1})$ where $\{\lambda\in\Lambda_1:h_\lambda\neq 1\}$ is well-ordered. We will call $H\ \wr_\omega \Lambda_1$ the \emph{lexicographic wreath product} of $H$ and $\Lambda_1$. Note that the restricted wreath product
is contained in the lexicographic wreath product.

\begin{theorem}
Let $H$ be a group, $\Lambda_0$ and $\Lambda_1$  \oag s, and $\theta:H\to\Aut^+(\Lambda_0)$ a homomorphism. Fix a $\theta$-affine action of $H$ on $\Lambda_0$ that preserves the orientation. Let
$\Lambda=\Lambda_0\times\mathcal{L}_{\lambda_1\in\Lambda_1}\Lambda_0$
and set $G=H\ {\wr_\omega}\ \Lambda_1$.
Then $\Lambda$ is an \oag. Define $\alpha:G\to\Aut^+(\Lambda)$ via $$\alpha_{(\lambda^*,(h_\lambda)_{\lambda\in \Lambda_1})}(\lambda',(\mu_{\lambda})_{\lambda\in\Lambda_1})
=(\lambda',(\theta_{h_{\lambda+\lambda^*}}\mu_{\lambda+\lambda^*})_{\lambda\in \Lambda_1}).$$
The action of $G$ on $\Lambda$ defined by $(\lambda^*,(h_{\lambda})_{\lambda\in \Lambda_1})\cdot(\lambda',(\mu_{\lambda})_{\lambda\in\Lambda_1})
=(\lambda'+\lambda^*,(h_{\lambda+\lambda^*}\mu_{\lambda+\lambda^*})_{\lambda\in H})$ is $\alpha$-affine. If the given action of $H$ is 
respectively \be\item free
\item rigid
\item 
essentially free
\ee 
then so is that of $G$.
\end{theorem}

\pf Recall that multiplication of elements of $G$ is performed as follows.
$$(\lambda^+,(k_\lambda))(\lambda^*,(h_\lambda))=(\lambda^++\lambda^*,(k_{\lambda-\lambda^*}h_\lambda))$$ 

One can check that $\left[(\lambda^+,(k_\lambda))(\lambda^*,(h_\lambda))\right]\cdot (\lambda',(\mu_\lambda))$ and $(\lambda^+,(k_\lambda))\left[(\lambda^*,(h_\lambda))\cdot (\lambda',(\mu_\lambda))\right]$ are both equal to $(\lambda^++\lambda^*+\lambda',(k_{\lambda+\lambda^+}h_{\lambda+\lambda^++\lambda^*}\mu_{\lambda+\lambda^++\lambda^*}))$. A similar calculation shows that $\alpha$ is a homomorphism.

To see that the action of $G$ on $\Lambda$ is $\alpha$-affine take $\bar{h}=(\lambda^*,(h_\lambda))\in G$ and $\bar{\lambda}'=(\lambda',(\mu'_\lambda))$ and $\bar{\lambda}''=(\lambda'',(\mu''_\lambda))\in\Lambda$ with $\bar{\lambda}'\geq\bar{\lambda}''$. Observe that $d_\Lambda(\bar{\lambda}',\bar{\lambda}'')=\bar{\lambda}'-\bar{\lambda}''$ while
\begin{eqnarray*}d_\Lambda(\bar{h}\bar{\lambda}' , \bar{h}\bar{\lambda}'')&=&(\lambda^*,(h_\lambda))\cdot(\lambda',(\mu'_\lambda))
-(\lambda^*,(h_\lambda))\cdot(\lambda'',(\mu''_\lambda))\\
&=&(\lambda'+\lambda^*,(h_{\lambda+\lambda^*}\mu'_{\lambda+\lambda^*}))
-(\lambda''+\lambda^*,(h_{\lambda+\lambda^*}\mu''_{\lambda+\lambda^*}))\\
&=& (\lambda'-\lambda'',\theta_{h_{\lambda+\lambda^*}}(\mu'_{\lambda+\lambda^*}-\mu''_{\lambda+\lambda^*}))\\ &=&\alpha_{\bar{h}}d_\Lambda(\bar{\lambda}',\bar{\lambda}'').
\end{eqnarray*}

If $\bar{h}\bar{\lambda}'=\bar{\lambda}'$ then $\lambda^*=0$ giving, for each $\lambda\in\Lambda_1$, $\mu_\lambda=h_{\lambda+\lambda^*}\mu_{\lambda+\lambda^*}=h_\lambda\mu_\lambda$. If the given action of $H$ on $\Lambda_0$ is free this forces $h_\lambda=1$ for all $\lambda$, whence $\bar{h}=1$.

Suppose now
that the given action of $H$ is rigid and $\bar{h}\bar{\lambda}'>\bar{\lambda}'$. To establish rigidity of the action of $G$ it suffices to show that $\bar{h}\bar{\lambda}''>\bar{\lambda}''$. If $\lambda^*>0$ this claim is obvious, so suppose that $\lambda^*=0$. Then taking the least $\lambda_0$ for which $h_{\lambda_0}$ does not fix $\Lambda_0$ pointwise
we have $h_{\lambda_0}\mu_{\lambda_0}>\mu_{\lambda_0}$ for all $\mu_{\lambda_0}\in\Lambda_0$, since the action of $H$ on $\Lambda$ is assumed rigid. It follows that $(0,(h_\lambda))\cdot(\lambda'',(\mu''_\lambda))=(\lambda'',(h_\lambda\mu''_\lambda))>(\lambda'',(\mu''_\lambda))$, as required. This establishes (i) and (ii), and (iii) is now immediate.\qed

Taking the free isometric (and hence essentially free) action of an \oag\ on itself and using induction, one deduces the following.

\begin{corollary}\label{iterated-wreath}
Let $\Lambda_i$ ($1\leq i\leq n$) be \oag s. Then the iterated wreath product $\Lambda_1\ \wr_\omega\ \Lambda_2\ \wr_\omega\ \cdots\ \wr_\omega\ \Lambda_n$ admits an essentially free affine action on $\Lambda$, viewed as a linear $\Lambda$-tree for some \oag\ $\Lambda$.
\end{corollary}

By a theorem of {\v{S}}mel{\cprime}kin \cite{Shmelkin} a free group in a product variety $\mathfrak{V}\mathfrak{U}$ embeds in the $\mathfrak{V}$-verbal wreath product of the free groups in the respective varieties. Note that in case $\mathfrak{V}$ is the variety of abelian groups, the $\mathfrak{V}$-verbal wreath product of two groups coincides with the restricted wreath product, which embeds in the lexicographic wreath product.
We therefore deduce the following result from Corollary~\ref{iterated-wreath}.

\begin{corollary}
Every free soluble group (of given derived length) admits an essentially free affine action on $\Lambda$ for some \oag\ $\Lambda$.
\end{corollary}

\section*{Acknowledgements}
I would like to thank Karel Dekimpe for highlighting the possible role of affine structures in the problems considered here,
and Alexei Miasnikov
and Martin Newell
for helpful conversations. I would also like to thank Florian Eisele and Yves de Cornulier for answering a question of mine on MathOverflow in relation to section 2.

\vspace{0.5cm}

\begin{minipage}[t]{3 in}
\noindent Shane O Rourke\\ Department of Mathematics\\
Cork Institute of Technology\\
Rossa Avenue\\ Cork\\ IRELAND
\\ \verb"shane.orourke@cit.ie"
\end{minipage}

\end{document}